\documentclass[reqno,a4paper,final]{amsart}
\newif\ifPDF                                                                  %
\ifx\pdfoutput\undefined\PDFfalse                                             %
\else\ifnum\pdfoutput >0\PDFtrue                                              %
     \else\PDFfalse                                                           %
     \fi                                                                      %
\fi                                                                           %
\usepackage{amssymb,amstext,amsthm,eucal,amscd}
\newtheorem{thm}{Theorem}

\newtheorem{prop}[thm]{Proposition}

\begin{document}
\title[Lorentzian spaces with a $\mathcal{T}_{1}\oplus\mathcal{T}_{3}$ structure]{Lorentzian homogeneous spaces admitting a homogeneous structure of type $\mathcal{T}_{1}\oplus\mathcal{T}_{3}$}
\author[P. Meessen]{Patrick Meessen}
\address[PM]{Physics Department, Theory, CERN, Geneva, Switzerland}
\email{Patrick.Meessen@cern.ch}
\date{\today}
\thanks{CERN-PH-TH/2005-006}
\begin{abstract}
 We show that a Lorentzian homogeneous space admitting 
 a homogeneous structure of type $\mathcal{T}_{1}\oplus\mathcal{T}_{3}$
 is either a (locally) symmetric space or a singular homogeneous plane
 wave. 
\end{abstract}
\maketitle
A theorem by Ambrose and Singer \cite{art:ambrose}, generalized to
arbitrary signature in \cite{art:oubinna}, states that on a reductive
homogeneous space, there exists a metric connection
$\overline{\nabla}=\nabla -S$, with $\nabla$ the Levi-Civit\`a
connection, that parallelizes the Riemann tensor $R$, and the
$(1,2)$-tensor $S$,  {\em i.e.\/} $\overline{\nabla}g
=\overline{\nabla}R =\overline{\nabla}S =0$.  Since a $(1,2)$-tensor
in $D\geq 3$ decomposes into 3 irreps of $\mathfrak{so}(D)$, one can
classify the reductive homogeneous spaces by the occurrence of one of
these irreps in the tensor $S$ \cite{boek:tricerri,art:gadea}. This
leads to 8 different classes, which range from the maximal, denoted by
$\mathcal{T}_{1}\oplus\mathcal{T}_{2}\oplus\mathcal{T}_{3}$, to the
minimal $ \{ 0 \}$. Clearly, homogeneous spaces of
type $ \{ 0 \}$ are just symmetric spaces. Moreover, also the
homogeneous spaces admitting a homogeneous structure of type
$\mathcal{T}_{i}$ ($i=1,2$ or $3$) have been characterized. For the
case at hand it is worth knowing that the homogeneous spaces with 
a $\mathcal{T}_{3}$ structure, for which $S$ corresponds to a 3-form,
are naturally reductive spaces \cite{boek:tricerri,art:gadea} and that
strictly Riemannian homogeneous $\mathcal{T}_{1}$ spaces are locally symmetric spaces \cite{boek:tricerri}.
Since a homogeneous structure of type $\mathcal{T}_{1}$ is defined
by an invariant vector field $\xi$, one must distinguish between
two cases in the Lorentzian setting: the non-degenerate case, for which 
$\xi$ is a space- or time-like vector, and the degenerate case, when
$\xi$ is a null vector. In the former case, Gadea and Oubi\~{n}a \cite{art:gadea}
proven that, analogously to the strictly Riemannian case, the space is locally symmetric.
In the degenerate case, A.~Montesinos~Amilibia \cite{art:amilibia} showed
that a homogeneous Lorentzian space admitting a degenerate $\mathcal{T}_{1}$ 
structure is a time-independent singular homogeneous plane wave
\cite{Blau:2002js}. A small calculation shows that a generic, {\em i.e.\/} time-dependent,
singular homogeneous plane wave admits a degenerate
$\mathcal{T}_{1}\oplus\mathcal{T}_{3}$ structure, see {\em e.g.\/} Appendix \ref{sec:SingHomWave}.
(By a (non-)degenerate $\mathcal{T}_{1}\oplus\mathcal{T}_{3}$ structure, we mean that the 
vector field $\xi$ characterizing the $\mathcal{T}_{1}$ contribution has (non-)vanishing norm.) 
This then automatically leads to the
question of whether the singular homogeneous plane waves exhaust the degenerate case in the 
$\mathcal{T}_{1}\oplus\mathcal{T}_{3}$ class. 
As we will see, the answer is affirmative.
\par
In the $\mathcal{T}_{1}\oplus\mathcal{T}_{3}$ case the homogeneous
structure is given by \cite{boek:tricerri}
\begin{displaymath}
  \overline{\nabla}_{X}Y \ -\ \nabla_{X}Y \ =\ - S_{X}Y \ =\ -T_{X}Y \
  -\ g (X,Y )\xi \ +\ \alpha ( Y ) X \; ,
\end{displaymath}
where we have defined $\alpha  ( X ) = g ( \xi ,X )$,
and $T_{X}Y (= -T_{Y}X)$ is the $\mathcal{T}_{3}$
contribution.  Since the metric and $S$ are parallel under
$\overline{\nabla}$, and $\xi$ is the contraction of $S$, it follows
that $\overline{\nabla}\xi = 0$ or, written differently:
\begin{displaymath}
  \nabla_{X}\xi \; =\; T_{X}\xi \ +\ \alpha ( X)\xi \ -\ \alpha (\xi ) X\; .
\end{displaymath}
This equation, together with the fact that $T$ is a 3-form, implies
that $\nabla_{\xi}\xi =0$, {\em i.e.\/}  $\xi$ is a geodesic vector.
\par
Given an isometry algebra $\mathfrak{g}$ ({\em i.e.\/} the Lie algebra of a Lie group acting 
transitively by isometries on a given homogeneous space), with a reductive split
$\mathfrak{g}=\mathfrak{m}+\mathfrak{h}$, where $\mathfrak{h}\subseteq
\mathfrak{so}(1,n+1)$ is the isotropy subalgebra, it is possible, and
usually done, to identify $\mathfrak{m}$ with $\mathbb{R}^{1,n+1}$;
the action of $\mathfrak{h}$ on $\mathfrak{m}$ can then be given by the
vector representation of $\mathfrak{so}(1,n+1)$ \cite{boek:kobayashi}.
This identification enables one to express the algebra in terms of $S$
and the curvature $\overline{R}$ as, limiting ourselves to the
$\mathfrak{m}\times\mathfrak{m}$ commutator,
\begin{equation}
  \label{eq:MMalg}
   [ X,Y ] \; =\; S_{X}Y \ -\ S_{Y}X \ +\ \overline{R} (
  X,Y ) \; ,
\end{equation}
where $S$ and $\overline{R}$ are evaluated at some point $p$. In the
above formula, $\overline{R}$ signals the presence of $\mathfrak{h}$
in $ [\mathfrak{m},\mathfrak{m} ]$. From now on, we only
consider this Lie algebra and all the relevant tensor fields are
evaluated at a specific point, even though this is not stated explicitly.
\par
Up to this point not too much has been said about $\mathfrak{h}$, and in
fact not too much can be said. It is known, however \cite{boek:kobayashi},
that a tensor field parallelized by $\overline{\nabla}$, when evaluated at a point corresponds to
an $\mathfrak{h}$-invariant tensor. Since in this
article we take $\xi$ (an $\mathfrak{h}$-invariant vector field as $\overline{\nabla}\xi =0$) 
to be non-vanishing, this means that $\mathfrak{h}\subseteq\mathfrak{so}(n+1)$
when $\xi$ is light-like, $\mathfrak{h}\subseteq\mathfrak{so}(1,n)$ when $\xi$ is space-like,
and $\mathfrak{h}\subseteq\mathfrak{iso}(n)$ when $\xi$ is null.
\par
Let us briefly outline the manner in which we arrive at our results: given a reductive
homogeneous space with reductive split $\mathfrak{g}=\mathfrak{m}+\mathfrak{h}$, the subalgebra
$\mathfrak{g}^{\prime}=\mathfrak{m}+ [ \mathfrak{m},\mathfrak{m} ]=\mathfrak{m}+\mathfrak{h}^{\prime}$ 
is an ideal of $\mathfrak{g}$.
It is this ideal, which is the Lie algebra of a Lie group still acting transitively, that we will consider;
we will say that an element
of $\mathfrak{h}$ appears in the algebra if it is an element of $\mathfrak{h}^{\prime}$. Given the
homogeneous structure, we can then, following Eq.\/ (\ref{eq:MMalg}), write down the maximal form of the
algebra compatible with the homogeneous structure. Since we are dealing with a Lie algebra, we can then
use the Jacobi identities to constrain the structure constants; after a redefinition of some generators
in $\mathfrak{m}$, corresponding to the choice of a different reductive split, this leads to a recognizable result.
Since the non-degenerate case is far less involved than the degenerate case, and gives a better idea of the 
manipulations used, it will be discussed before the degenerate case.
\section{The non-degenerate case}
Let $\mathfrak{m}$ be spanned by the generators $V$ and $Z_{i}$ ($i=1,\ldots ,n$), which in this case we can take to satisfy
\begin{displaymath}
  \begin{array}{lclclcl}
      \langle V,V\rangle & =& \aleph   &\hspace{.3cm},\hspace{.3cm}& \alpha ( V ) & =& \lambda\ =\ \aleph |\lambda | \; ,\\
      \langle Z_{i},Z_{j}\rangle & =& \eta_{ij} & ,& \alpha ( Z_{i} ) & =& 0 \; , \\
      \langle V, Z_{i}\rangle & = & 0 & ,& & & 
  \end{array}
\end{displaymath}
where $\aleph =\pm 1$ distinguishes between the time-like (for $\aleph =-1$) and the space-like (for $\aleph =1$) cases 
and $\eta =\mathrm{diag}( -\aleph , 1,\ldots ,1)$.
As is mentioned in the introduction, $\mathfrak{h}$ is contained in either $\mathfrak{so}(n+1)$ (for $\aleph =-1$)
or $\mathfrak{so}(1,n)$ (for $\aleph =1$) and the relevant non-vanishing commutation relations are
\begin{displaymath}
  \begin{array}{lcl}
   [ M_{ij},M_{kl} ] & =& \eta_{jk}\ M_{il}
                                  -\eta_{ik}\ M_{jl}
                                  +\eta_{jl}\ M_{ki}
                                  -\eta_{il}\ M_{kj} \; ,\\ {}
   [ M_{ij},Z_{k} ] & =&  \eta_{jk}\ Z_{i}
                                  -\eta_{ik}\ Z_{j} \; .{}
  \end{array}
\end{displaymath}
Once again, let us stress that not every $M$ needs appear, but the elements of $\mathfrak{h}^{\prime}$ 
can be written as combinations of the $M$'s, and their commutation relations are induced by
the ones above.
\par
With respect to the chosen basis we can decompose $2 T_{V}Z_{i}= F_{i}{}^{j}Z_{j}$ and 
$2 T_{Z_{i}}Z_{j} = \aleph F_{ij}V \ +\ C_{ij}{}^{k}Z_{k}$, which allows
us to write
\begin{eqnarray}
   [ V,Z_{i} ] & =& \lambda\ Z_{i} \ +\ F_{i}{}^{j}Z_{j} \ +\ \overline{R} ( V,Z_{i} ) 
  \; , \nonumber \\{}
  & & \nonumber \\{}
   [ Z_{i},Z_{j} ] & =& \aleph F_{ij}V \ +\ C_{ij}{}^{k}Z_{k} \ +\ \overline{R} ( Z_{i},Z_{j} ) 
  \; . \nonumber {}
\end{eqnarray}
\par
Let us then, following the strategy outlined above, check the Jacobi identities. The first one 
is the $(V,Z_{i},Z_{j})$ identity, which leads to $F=0$ and
\begin{eqnarray}
  \label{eq:nondegenTZZidentity1}
  \textstyle{\lambda\over 2} C_{ijk} & =& R_{jik} \; -\; R_{ijk} \, ,\\
  \label{eq:nondegenTZZidentity2}
  2\lambda\ S_{ij}{}^{mn} & =& C_{ij}{}^{k}\ R_{k}{}^{mn} \; , 
\end{eqnarray}
where we expanded $\overline{R} ( V,Z_{i} ) = R_{i}{}^{mn}\ M_{mn}$ and 
$\overline{R} ( Z_{i},Z_{j} ) = S_{ij}{}^{mn}\ M_{mn}$.
Since $F=0$ we can redefine 
\begin{displaymath}
  Y_{i} \; =\; Z_{i} \; +\; \lambda^{-1}\ R_{i}{}^{mn}M_{mn} \; ,
\end{displaymath}
from which we trivially find
\begin{displaymath}
  \label{eq:nondegTY}
   [ V, Y_{i} ] \; =\; \lambda\ Y_{i} \; , 
\end{displaymath}
which at once implies that $C=0$, by Eq.~(\ref{eq:nondegenTZZidentity1}), and also that $S=0$ thanks to Eq.~(\ref{eq:nondegenTZZidentity2}).
So the, quite remarkable, result is that a Lorentzian homogeneous space admitting a non-degenerate homogeneous structure
of type $\mathcal{T}_{1}\oplus\mathcal{T}_{3}$, also admits a non-degenerate $\mathcal{T}_{1}$ structure.
Combining this with the results of Gadea and Oubi\~{n}a \cite{art:gadea}, we have proven the following result. 
\begin{prop}
\label{prop1}
A connected homogeneous Lorentzian
space admitting a non-degenerate $\mathcal{T}_{1}\oplus\mathcal{T}_{3}$ structure is a locally symmetric space.
\end{prop}
\section{The degenerate case}
In the degenerate case we can choose the generators $U$, $V$ and $Z_{i}$ ($i=1,\ldots , n$) spanning $\mathfrak{m}$ such
that $\alpha  ( U ) = \lambda\neq 0$, $\alpha ( V ) =\alpha ( Z_{i} ) =0$.
The invariant norm is then $\langle U,V\rangle =1$ and $\langle Z_{i},Z_{j}\rangle = \delta_{ij}$
and we decompose the $\mathcal{T}_{3}$ contribution to $S$ as
\begin{displaymath}
  \label{eq:DegenerateSplit}
  \begin{array}{lclclcl}
      2T ( U,V,Z_{i} ) & =& W_{i} &\hspace{.2cm},\hspace{.2cm}&
      2T ( U,Z_{i},Z_{j} ) & =& F_{ij} \; , \\
      & & & & & & \\
      2T ( Z_{i},Z_{j},Z_{k} ) & =& C_{ijk} &\hspace{.2cm},\hspace{.2cm}&
      2T ( V,Z_{i},Z_{j} ) & =& \aleph_{ij} \; , 
  \end{array}
\end{displaymath}
where $F$, $\aleph$ and $C$ are totally antisymmetric. Given these abbreviations we can write
the most general $\mathfrak{m}\times\mathfrak{m}$ commutators as
\begin{eqnarray}
  \label{eq:DegMostGen1}
   [ U,V ] & =& \lambda\ V \, +\, W^{i}\ Z_{i}\, +\, \overline{R} ( U,V ) \; , \nonumber \\
  & & \nonumber \\
  \label{eq:DegMostGen2}
   [ U,Z_{i} ] & =& \lambda\ Z_{i} \, +\, F_{i}{}^{j}\ Z_{j} \, -\, W_{i}\ U
                             \, +\, \overline{R} ( U,Z_{i} ) \; , \nonumber \\
  & & \nonumber \\
   \label{eq:DegMostGen3}
   [ V,Z_{i} ] & =& W^{i}\ V \, +\, \aleph_{i}{}^{j}\ Z_{j}\, +\, \overline{R} ( V,Z_{i} ) \; , \nonumber \\
  & & \nonumber \\
  \label{eq:DegMostGen4}
   [ Z_{i},Z_{j} ] & =& \aleph_{ij}\ U \, +\, F_{ij}\ V \, +\, C_{ijk}\ Z^{k}
                                 \, +\, \overline{R} ( Z_{i},Z_{j} ) \; ,\nonumber
\end{eqnarray}
where the various $\overline{R}$ need to be expanded in terms of the generators of $\mathfrak{h}$. Since $\xi$
is null, we see that $\mathfrak{h}\subseteq\mathfrak{iso}(n)$, which we take to be spanned by $\overline{Z}_{i}$
and $M_{ij}$ with commutation relations
\begin{eqnarray}
 \label{eq:IsoAction}
   [ M_{ij},M_{kl} ] & =& \delta_{jk}M_{il}
                                  -\delta_{ik}M_{jl}
                                  +\delta_{jl}M_{ki}
                                  -\delta_{il}M_{kj} \; , \nonumber \\ {}
   [ M_{ij},\overline{Z}_{k} ] & =&
                                  \delta_{jk}\overline{Z}_{i}
                                  -\delta_{ik}\overline{Z}_{j} \; ,\nonumber \\{}
   [ M_{ij},Z_{k} ] & =&
                                  \delta_{jk}Z_{i} -\delta_{ik}Z_{j}
                                  \; ,\nonumber \\ {}
   [ U,\overline{Z}_{i} ] & =& Z_{i} \; , \nonumber \\ {}
   [ Z_{i},\overline{Z}_{j} ] & =& -\delta_{ij}V \; , \nonumber {}
\end{eqnarray}
where it should be kept in mind that not all elements of
$\mathfrak{iso}(n)$ need appear.
\par
We can then once again start to recover the information contained in the Jacobi identities:
the $ ( U,V,Z )$ Jacobi identity reads
\begin{eqnarray}
  \label{eq:PreUVZ}
  0 & =& -2\lambda W_{i}\ V \, -\,  \{
                                    \lambda\aleph_{ij} +F_{i}{}^{k}\aleph_{kj} + F_{j}{}^{k}\aleph_{ik} + W^{k}C_{kij}  
                                   \}\ Z^{k} \nonumber \\
 & & \nonumber \\
    & & - [ \overline{R} ( U,V ),Z_{i} ] \ -\  [ \overline{R} ( V,Z_{i} ),U ] \nonumber \\
 & & \nonumber \\
 & & +\aleph_{i}{}^{j}\overline{R} ( U,Z_{i} ) -2\lambda\overline{R} ( V,Z_{i} )
        -F_{i}{}^{j}\overline{R} ( V,Z_{i} ) + W^{j}\overline{R} ( Z_{i},Z_{j} ) \; .
\end{eqnarray}
Cancellation of the $V$ contribution then means that $\overline{R} ( U,V ) = -2\lambda W^{i}\overline{Z}_{i} + Y^{ij}M_{ij}$,
which at once means that $W$ can only be non-zero for those directions for which a $\overline{Z}$ appears. Specifically,
should none appear, then $W=0$. Let us then split the index $i$ into some indices $a$ and $I$, such that 
the $\overline{Z}_{a}$ do appear whereas the $\overline{Z}_{I}$ do not. 
\par
Having made the split, we can investigate the implication of having the null-boosts in the algebra. Let
us start by looking at the $ ( U, Z_{i},\overline{Z}_{a} )$ Jacobi: a small calculation then
shows that this implies
\begin{eqnarray}
  \label{eq:PreUZbarZ}
  0 & =& -\aleph_{ia}\ U \ -\delta_{ia}W^{i}\ Z_{i} \ +\ W_{i}Z_{a} \ +\ C_{aik}\ Z^{k} \nonumber \\
 & & \nonumber \\
 & & - [ \overline{R} ( U,Z_{i} ), \overline{Z}_{a} ]
     \ -\ \delta_{ia}\overline{R} ( U,V )
     \ -\ \overline{R} ( Z_{i},Z_{a} ) \; . \nonumber
\end{eqnarray}
In order for the above to be true we must have that $\aleph_{ai}=C_{aij}=0$ and that $W$ can be non-zero
only if no or only one $\overline{Z}$ appears in $\mathfrak{h}$. As was said above, the no-case already
implies that $W=0$, so we had better have a look at the case of one appearing null boost. For this 
we are helped by the $\mathfrak{h}$-part of the above equation. Clearly in the case when we are dealing with
only one $\overline{Z}$, this amounts to the statement that 
$ [ \overline{R} ( U,Z_{a} ),\overline{Z}_{a} ] = -\overline{R} ( U,V )$, which,
since there is no rotation in $\mathfrak{so}(n)$ that can take $Z_{a}$ to $Z_{a}$, means
that $\overline{R} ( U,V )=0$, and hence that $W_{a}=0$. This then means that in all cases we
have $W=0$.
\par
Continuing with the analysis, one can see that the $ ( Z_{i},Z_{j},\overline{Z}_{a} )$ Jacobi
leads to 
\begin{eqnarray}
  \aleph_{ij}\ Z_{a} & =& \delta_{ja}\aleph_{i}{}^{k}Z_{k} \ -\ \delta_{ia}\aleph_{j}{}^{k}Z_{k}\; , \nonumber \\ {}
 & & \nonumber \\ {}
   [\overline{R} ( Z_{i},Z_{j} ),\overline{Z}_{a} ] & =& 
      \delta_{ja}\overline{R} ( U,Z_{i} )
     \ -\ \delta_{ia}\overline{R} ( U,Z_{j} ) \; . \nonumber {}
\end{eqnarray}
Then, using the fact that $\aleph_{ia}=0$, one then sees that $\aleph_{IJ}=0$ and that hence $\aleph_{ij}=0$
when $\mathfrak{h}$ includes some null boost. In the case when there is no $\overline{Z}$, the relevant
information can be obtained by picking out the $V$ component in the $ ( V,Z_{i},Z_{j} )$ Jacobi:
this implies that $\lambda \aleph_{ij} = F_{i}{}^{k}\aleph_{kj} + F_{j}{}^{k}\aleph_{ik}$, which after
contraction leads to $\lambda \aleph_{ij}\aleph^{ij}=0$ and thus implies that $\aleph =0$. 
\par
The $\mathfrak{h}$-part of Eq.~(\ref{eq:PreUVZ}) then implies that 
$2\lambda \overline{R} ( V,Z_{i} ) = - F_{i}{}^{j}\overline{R} ( V,Z_{j} )$, so that
$\overline{R} ( V,Z_{i} ) =0$. In order to then identically satisfy Eq.~(\ref{eq:PreUVZ})
we must have $ [ \overline{R} ( U,V ),Z_{i} ] =0$, so that $\overline{R} ( U,V ) =0$.
\par
Summarizing the results obtained thus far, we find that the non-trivial $\mathfrak{m}\times\mathfrak{m}$-commutators, 
scaling $U$ in such a way that $\lambda =1$ and decomposing the various $\overline{R}$'s, are
\begin{eqnarray}
  \label{eq:StartingPoint}
   [ U,V ] & =& V \; , \nonumber \\ {}
   [ U, Z_{i} ] & =&  ( F+\delta )_{ij}Z_{j} +
                              h_{ij}\overline{Z}_{j} +
                              \textstyle{1\over 2}R_{ijk}M_{jk} \; ,\nonumber \\ {}
   [ Z_{i},Z_{j} ] & =& F_{ij}V +
                              C_{ijk}Z_{k} + S_{ijk}\overline{Z}_{k} +
                              N_{ijkl}M_{kl} \; .\nonumber {}
\end{eqnarray}
Let us then continue our analysis of the Jacobi identities:
the $ ( U,Z_{i},Z_{j} )$ Jacobi implies
\begin{eqnarray}
  h_{ij} & =& A_{(ij)} \, -\, \textstyle{1\over 2}F_{ij}  \; , \nonumber \\ 
  C_{ijk}h_{kl} & =&  ( F+\delta )_{ik} S_{kjl} 
                \ +\  ( F+\delta )_{jk} S_{ikl} \; , \nonumber \\ 
  \textstyle{1\over 2}C_{ijk}R_{kmn} & =&  ( F+\delta )_{ik} N_{kjmn} 
                                     \ +\  ( F+\delta )_{jk} N_{ikmn} \label{eq:UZZ-3}\; ,\\ 
  S_{ijk} +  R_{ijk} -R_{jik} & =& \delta_{F}C_{ijk} \ +\ C_{ijk}  \label{eq:UZZ-4}\; ,
\end{eqnarray}
where we defined
\begin{displaymath}
  \delta_{F}C_{ijk} \; =\; F_{il}C_{ljk} \ +\ F_{jl}C_{ilk} \ +\
  F_{kl}C_{ijl} \; .
\end{displaymath}
From Eq.~(\ref{eq:UZZ-4}) one sees that $S$ must be totally
antisymmetric.  Denoting by $\mathfrak{S}_{(ijk)}$ the sum over the permutations $(ijk)$, $(jki)$ and $(kij)$,  
the $ ( Z_{i},Z_{j},Z_{k} )$ Jacobi results in
\begin{eqnarray}
  0 & =& \mathfrak{S}_{(ijk)}\ C_{jkl}S_{ilm} \label{eq:ZZZ-1}\; , \nonumber \\
  0 & =& \mathfrak{S}_{(ijk)}\ C_{jkl}N_{ilmn} \label{eq:ZZZ-2}\; , \nonumber \\
  0 & =& \mathfrak{S}_{(ijk)} [ C_{jkl}C_{ilm} \ +\ 2N_{jkim}
   ] \label{eq:ZZZ-3}\; , \nonumber
\end{eqnarray}
and also, since $S$ is totally antisymmetric, 
\begin{equation}
  \label{eq:ZZZ-4}
  3 S \; =\; \delta_{F}C \; .
\end{equation}
Of course, if a $\overline{Z}_{a}$ occurs in $ [
\mathfrak{m},\mathfrak{m} ]$, then the $ (
U,Z_{i},\overline{Z}_{a} )$ Jacobi implies that
\begin{eqnarray}
  \label{eq:UZbarZ-1}
  C_{iaj} & =& 0 \; , \nonumber \\ 
  S_{iaj} & =& R_{iaj} \label{eq:UZbarZ-2}\; ,\\
  N_{iakl} & =& 0 \label{eq:UZbarZ-3}\; .
\end{eqnarray}
\par
Let us then, as before, split the indices $i$ into $(a,I)$, where the $\overline{Z}_{a}$'s occur
but the $\overline{Z}_{I}$'s do not. This means by assumption that $h_{iI}=0$, which
implies $2A_{aI}=F_{aI}$, $A_{IJ}=0=F_{IJ}$ and $S_{ijI}=0$,
which implies that only $S_{abc}$ is non-zero. Furthermore, we then
see that only $C_{IJK}$ is non-vanishing.  Together with
Eq.~(\ref{eq:ZZZ-4}), this then implies that $S=0$, and we get the
extra constraint
\begin{equation}
  \label{eq:Restr-1}
    F_{aI}C_{IJK} \; =\; 0 \; .
\end{equation}
This last constraint also follows from the $ (
Z_{i},Z_{j},\overline{Z}_{a} )$ Jacobi, which also tells us that
$N_{ijal}=0$.
\par
Eq.~(\ref{eq:UZbarZ-2}) then implies that only $R_{IJK}$ and $R_{aJK}$
are non-vanishing, and from Eq.~(\ref{eq:UZbarZ-3}) we find that only
$N_{IJmn}$ can be non-zero. We can calculate $R_{aJK}$ from
Eq.~(\ref{eq:UZZ-4}), which then gives $R_{aIJ} = F_{aK}C_{KIJ} = 0$
because of Eq.~(\ref{eq:Restr-1}).  The same equation then states
$R_{IJK}-R_{JIK}=C_{IJK}$, which by means of Eq.~(\ref{eq:UZZ-3}) then
also implies that only the $N_{IJKL}$ can be non-vanishing.
\par
Let us define the generator
\begin{displaymath}
  \label{eq:ZIRedef}
  Y_{I} \; =\; Z_{I} \ -\ F_{Ia}\overline{Z}_{a} \; ,
\end{displaymath}
from which we can then derive that the algebra takes on the form
\begin{eqnarray}
  \label{eq:NewAlgebra}
   [ U, Z_{a} ] & =&  ( F+\delta )_{ab}Z_{b} \ +\
                             ( A_{ab} -\textstyle{1\over
                            2}F_{ab} ) \overline{Z}_{b} \; ,\nonumber \\ {}
   [ Z_{a},Z_{b} ] & =& F_{ab}V\; , \nonumber  \\ {}
   [ U, Y_{I} ] & =& Y_{I} \ +\ \textstyle{1\over 2}R_{IJK}M_{JK} \; , \nonumber \\ {}
   [ Y_{I},Y_{J} ] & =& C_{IJK}Y_{K} \ +\ N_{IJKL}M_{KL} \; ,\nonumber {}
\end{eqnarray}
so that the $a$- and the $I$-sectors decouple from each other.
\par
Restricting ourselves to the $I$-sector and further defining
\begin{displaymath}
  \label{eq:YIRedef}
  W_{I} \; =\; Y_{I} \ +\ \textstyle{1\over 2}R_{IJK}M_{JK} \; ,
\end{displaymath}
we immediately find $ [ U,W_{I} ] \ =\ W_{I}$;
calculating the remaining commutator, we find
\begin{displaymath}
   [ W_{I},W_{J} ] \ =\  ( C_{IJK} - R_{IJK} +
         R_{JIK} ) Y_{K} \ + \ldots\ ,
\end{displaymath}
where the $\ldots$ stands for terms in $M_{JK}$. Using now
Eq.~(\ref{eq:UZZ-4}), we see that this redefinition trivializes $C$,
and by way of Eq.~(\ref{eq:UZZ-3}), also $N$.
\par
At this point, the only difference between the algebra we deduced and the 
generic singular homogeneous plane wave algebra in Eq.~(\ref{eq:sing_hom_waves_algebra})
are the null boosts in the $I$-sector, that is a generator one would call $\overline{W}_{I}$.
It is, however, always possible to extend our algebra to an algebra that does contain them; 
in fact this follows immediately from the consistency of the singular homogeneous plane wave algebra.
Putting everything together, one sees that we obtain the isometry algebra
of a generic singular homogeneous plane wave in Eq.~(\ref{eq:sing_hom_waves_algebra})
by, basically, choosing a different reductive split of the same algebra.
Thus we have proven the next theorem.
\begin{thm}
The underlying geometry of a connected homogeneous Lorentzian space
that admits a degenerate $\mathcal{T}_{1}\oplus\mathcal{T}_{3}$ structure
is that of a singular homogeneous plane wave.
\end{thm}

\par\noindent
{\bf Note added}: The author recently became aware of \cite{art:pastore}, where Proposition \ref{prop1}
is proven for the Riemannian case. The reasonings leading to Proposition \ref{prop1}, however,
only depend on the non-degeneracy of the $\mathcal{T}_{1}$ contribution and not on the signature
of the space. This means that Proposition \ref{prop1} also holds in the pseudo-Riemannian case.  
\section*{Acknowledgments}
The author would like to thank J.~Figueroa-O'Farrill, R.~Hern\'andez and
S.~Philip for very useful discussions, and S. Vascotto for improving the readability
of the text.
\appendix{
\section{Singular homogeneous plane waves}
\label{sec:SingHomWave}
A global coordinate system for the singular homogeneous plane waves is
defined by the data\footnote{ This form of the metric is related to
the one in \cite[Eq.~(2.51)]{Blau:2002js} by the transformations
$x^{+}= e^{-z}$, $x^{-}= -e^{z}s$, $\vec{z}=\vec{x}$, $A_{0}= 2H$ and $f = -F$.
}
\begin{displaymath}
  \label{eq:sing_hom_wave}
  \begin{array}{lcl}
    e^{+} & =& dz \; ,\\
    & &  \\
    e^{-} & =& ds +   [ \vec{x}^{T}e^{zF}He^{-zF}\vec{x} +s ] dz \; ,\\
    & & \\ 
    e^{i} & =& dx^{i} \; ,
  \end{array}
\end{displaymath}
where the metric is defined by $\eta_{+-} =1$ and $\eta_{ij} = \delta_{ij}$.
This class of metrics admits a homogeneous structure given by the components
\begin{displaymath}
  \label{eq:sing_hom_wave_hom_struct}
  S_{++-} \ =\ -1 \hspace{.3cm},\hspace{.3cm}
  S_{+ij} \ =\ F_{ij} \hspace{.3cm},\hspace{.3cm}
  S_{i+j} \ =\ -\delta_{ij} - F_{ij} \; ,
\end{displaymath}
which corresponds to a degenerate $\mathcal{T}_{1}\oplus\mathcal{T}_{3}$ structure.  
\par
The isometry algebra, apart from possible rotations that appear as automorphisms of the algebra,
can be found to be \cite{Blau:2002js}
\begin{equation}
  \label{eq:sing_hom_waves_algebra}
  \begin{array}{lclclcl}
      [ U, V ] & =& V &\hspace{.3cm},\hspace{.3cm}& {}
      [ \overline{X}_{i},\overline{X}_{j} ] & =& 0 \; \\ {}
     & & & & & & \\ {}
      [ X_{i}, X_{j} ] & =& 2F_{ij}\ V  {}   & ,& {}
      [ X_{i}, \overline{X}_{j} ] & =& -\delta_{ij}\ V \; \\ {}
     & & & & & & \\{}
      [ U, \overline{X}_{i}  ] & =& X_{i} & ,& {}
      [ U, X_{i}  ] & =&  [ 2H -F ]_{ij}\overline{X}_{j} 
                                  \ +\ 
                                  [ \delta + 2F ]_{ij}
                                  X_{j} \; . {}
  \end{array}
\end{equation}
}


\begin{thebibliography}{99}
%
\bibitem{art:ambrose} W.~Ambrose, I.~M.~Singer:
``{\em On homogeneous Riemannian manifolds}'',
Duke Math. J. {\bf 25} (1958), 647--669.
%
\bibitem{art:oubinna} P.~M.~Gadea, J.~A.~Oubi\~{n}a:
``{\em Homogeneous pseudo-Riemannian structures and homogeneous almost para-Hermitian structures}'',
Houston J. Math. {\bf 18} (1992), 449--465.
%
\bibitem{boek:tricerri} F.~Tricerri, L.~Vanhecke:
``{\em Homogeneous structures on Riemannian manifolds}'',
London Math. Soc. Lecture Note Ser. {\bf 83} (1983), 1--125.
%
\bibitem{art:gadea} P.~M.~Gadea, J.~A.~Oubi\~{n}a:
``{\em Reductive homogeneous pseudo-Riemannian manifolds}'',
Monatsh. Math. {\bf 124} (1997), 17--34.
%
\bibitem{art:amilibia} A.~Montesinos~Amilibia:
``{\em Degenerate homogeneous structures of type $\mathcal{S}_{1}$ on pseudo-Riemannian manifolds}'',
Rocky Mountain J. Math. {\bf 31} (2001), 561--579.
%
\bibitem{Blau:2002js}
M.~Blau, M.~O'Loughlin:
``{\em Homogeneous plane waves}'',
Nuclear Phys. B {\bf 654} (2003), 135--176.
%
\bibitem{boek:kobayashi} S.~Kobayashi, K.~Nomizu,
``{\em Foundations of differential geometry}'',
Wiley (1963 and 1969).
%
\bibitem{art:pastore} A.~M.~Pastore,
``{\em On the homogeneous Riemannian structure of type $\mathcal{T}_{1}\oplus\mathcal{T}_{3}$}'',
Geom. Dedicata {\bf 30}(1989), 235--246.
%
\end{thebibliography}
\end{document}